\documentclass[12pt,oneside]{amsart}
\usepackage{amssymb,verbatim,enumerate,amsmath,ifthen,cite}
\usepackage[mathscr]{eucal}
\usepackage[utf8]{inputenc}
\usepackage[T1]{fontenc}
\usepackage[marginparwidth=2.4cm]{geometry}
\newtheorem{theorem}{Theorem}[section]
\newtheorem*{theorem*}{Theorem}
\def\Thm#1#2{\ifthenelse{\equal{#1}{*}}{\begin{theorem*}#2\end{theorem*}}
             {\begin{theorem}\label{T#1}#2\end{theorem}}}
\newtheorem{Atheorem}{Theorem}

\def\THM#1#2{\begin{Atheorem}\label{T#1}#2\end{Atheorem}}
\def\thm#1{Theorem~\ref{T#1}}

\newtheorem{proposition}[theorem]{Proposition}
\newtheorem*{proposition*}{Proposition}
\def\Prp#1#2{\ifthenelse{\equal{#1}{*}}{\begin{proposition*}#2\end{proposition*}}
             {\begin{proposition}\label{P#1}#2\end{proposition}}}

\newtheorem{corollary}[theorem]{Corollary}
\newtheorem*{corollary*}{Corollary}
\def\Cor#1#2{\ifthenelse{\equal{#1}{*}}{\begin{corollary*}#2\end{corollary*}}
             {\begin{corollary}\label{C#1}#2\end{corollary}}}

\newtheorem{lemma}[theorem]{Lemma}
\newtheorem*{lemma*}{Lemma}
\def\Lem#1#2{\ifthenelse{\equal{#1}{*}}{\begin{lemma*}#2\end{lemma*}}
             {\begin{lemma}\label{L#1}#2\end{lemma}}}
\def\lem#1{Lemma~\ref{L#1}}
\newtheorem{Alemma}{Lemma}

\theoremstyle{definition}
\newtheorem{remark}[theorem]{Remark}
\newtheorem*{remark*}{Remark}
\def\Rem#1#2{\ifthenelse{\equal{#1}{*}}{\begin{remark*}\rm #2\end{remark*}}
             {\begin{remark}\label{R#1}\rm #2\end{remark}}}

\newtheorem{example}[theorem]{Example}
\newtheorem*{example*}{Example}
\def\Exa#1#2{\ifthenelse{\equal{#1}{*}}{\begin{example*}\rm #2\end{example*}}
             {\begin{example}\label{Ex#1}\rm #2\end{example}}}

\numberwithin{equation}{section}
\def\eq#1{{\rm(\ref{#1})}}
\def\Eq#1#2{\ifthenelse{\equal{#1}{*}}
  {\begin{equation*}\begin{aligned}[]#2\end{aligned}\end{equation*}}
  {\begin{equation}\begin{aligned}[]\label{#1}#2\end{aligned}\end{equation}}}

\numberwithin{equation}{section}
\def\eq#1{{\rm(\ref{#1})}}
\def\Eq#1#2{\ifthenelse{\equal{#1}{*}}
  {\begin{equation*}\begin{aligned}[]#2\end{aligned}\end{equation*}}
  {\begin{equation}\begin{aligned}[]\label{#1}#2\end{aligned}\end{equation}}}

\def\A{\mathscr{A}}
\def\M{\mathscr{M}}

\def\P{\mathscr{P}}

\def\C{\mathcal{C}}

\newcommand\R{\mathbb{R}}
\newcommand\N{\mathbb{N}}
\newcommand\Cyc{\textrm{Cyc}}
\newcommand{\QA}[1]{\A_{#1}}

\begin{document}

\title[On subadditive quasi-arithmetic means]{On subadditive quasi-arithmetic means}

\author[Zsolt P\'ales]{Zsolt P\'ales}
\address{%
Institute of Mathematics, 
University of Debrecen,
Pf.\ 400, 4002 Debrecen, Hungary}
\email{pales@science.unideb.hu}
\author{Pawe{\l} Pasteczka}
\address{Institute of Mathematics, University of Rzesz\'ow, Pigonia 1, 35-310 Rzesz\'ow, Poland}
\email{ppasteczka@ur.edu.pl}


\begin{abstract}
Let $f\colon \R_+\to\R$ be a continuous and strictly monotone function. In the main result of this paper, we show that, for a fixed $n\geq 2$, the $n$-variable mean 
$\mathscr{A}_f \colon \R_+^n \to \R_+$ defined by
$$
\mathscr{A}_f(x_1,\dots,x_n):=f^{-1} \bigg( \frac{f(x_1)+\cdots+f(x_n)}n \bigg)
$$
is subadditive if and only if $f$ is differentiable with a continuously semi-differentiable and nonvanishing first derivative, and there exists an $\alpha\in[0,\infty]$ such that $f''_+:=(f')'_+$ is positive on $(0,\alpha)$ and $f''_+=0$ on $[\alpha,\infty)$, furthermore, $\frac{f'}{f''_+}$ is increasing and superadditive on $(0,\alpha)$.
\end{abstract}

\maketitle

\section{Introduction}

We are going to recall first several basic facts about quasi-arithmetic means.
Given a continuous and strictly monotone function $f \colon I \to \R$, the \emph{quasi-arithmetic mean} $\QA{f} \colon \bigcup_{n=1}^\infty I^n \to I$ is defined by
\Eq{*}{
\QA{f}(x):=f^{-1} \bigg( \frac{f(x_1)+\cdots+f(x_n)}n \bigg),\qquad\text{where }n \in \N,\,x=(x_1,\dots,x_n) \in I^n.
}
Knopp \cite{Kno28} observed that, for $I=\R_+$ and the functions $\pi_p(x):=x^p$ ($p\ne 0$) and $\pi_0(x):=\ln x$, the quasi-arithmetic mean $\QA{\pi_p}$ coincides with the \emph{$p$-th power mean} $\P_p$. The class of quasi-arithmetic means was axiomatized in the 1930s \cite{Kol30,Nag30,Def31}. It is also well-known that the homogeneous quasi-arithmetic means are exactly the power means.

A particularly important subclass of this family (containing, in particular, the power means) consists of means generated by $\C^2$ functions whose first derivative never vanishes. By applying Jensen’s inequality, it is straightforward to verify that, for $f,g$ possessing this regularity, the comparison inequality $\QA{f}\le\QA{g}$ is equivalent to
\[
\tfrac{f''}{f'}\le\tfrac{g''}{g'},
\]
i.e., to the pointwise comparability of two single-variable functions. For a detailed treatment of this connection, see Mikusi\'nski \cite{Mik48}. The operator $f\mapsto\tfrac{f''}{f'}$ has also been employed in several contexts by Pasteczka \cite{Pas13,Pas15d,Pas18c,Pas20a}.

Minkowski- and Hölder-type inequalities have been considered in the context of quasi-arithmetic means by several authors, cf.\ Beck \cite{Bec70}, Beckenbach--Bellmann \cite{BecBel61}, Gr\"unwald--P\'ales \cite{GruPal24}, Hardy--Littlewood--Pólya \cite{HarLitPol34}, Losonczi \cite{Los83a}, P\'ales \cite{Pal86a, Pal99a}. Recently, we have extended the notion of quasi-arithmetic means \cite{PalPas26} and in a subsequent paper \cite{PalPas2511}, we characterized the validity of Minkowski- and Hölder-type inequalities in this more general context. 

The problem of convexity of quasi-arithmetic means was first considered by the authors in \cite{PalPas18a} to obtain sufficient conditions for Hardy-type inequalities for quasi-arithmetic means. Then a comprehensive study of convexity and concavity of weighted quasi-arithmetic means was given by Chudziak--Głazowska--Jarczyk--Jarczyk \cite{ChuGlaJarJar19}. Finally, the results obtained by the authors in \cite{PalPas21a} established the following ultimate characterization of the convexity of quasi-arithmetic means. 

\THM{*}{Let $n\geq2$ be fixed and $f \colon I \to \R$ be continuous and strictly monotone. Then the $n$-variable mean $\A_f|_{I^n}$ is Jensen convex if and only if $f$ is twice continuously differentiable with a nonvanishing first derivative and either $f''$ is identically zero on $I$, or $f''$ is nowhere zero and $f'/f''$ is positive and convex on $I$.}

The convexity of more general means, such as Bajraktarević means and quasideviation means was characterized in \cite{PalPas23b}.

A surprising feature of this result is that convexity of a quasi-arithmetic mean implies very strong regularity properties on its generator. Furthermore, the convexity of the mean is equivalent to the convexity of a certain single variable function.

Motivated by this result, the main aim of this paper is to show that the subadditivity of a quasi-arithmetic mean also yields strong regularity properties of the generator function. Moreover, we show that such a property is equivalent to the superadditivity of a certain single variable function.

\section{Semi-differentiability, Monotonicity and Convexity}

Let $(X,\|\cdot\|)$ and $(Y,\|\cdot\|)$ be real normed spaces and let $D\subseteq X$ be nonempty and open.
Given a function $f:D\to Y$, we say that $f$ is \emph{semi-differentiable at $x\in D$} if, for all $u\in X$, the one-sided limit, called the \emph{directional derivative of $f$ at $x$ in direction $u$},
\Eq{*}{
  f'(x;u):=\lim_{t\to0+}\frac{f(x+tu)-f(x)}{t}. 
}
exists. If $f$ is semi-differentiable at $x$ and the map $u\mapsto f'(x;u)$ is linear, then we say that $f$ is \emph{differentiable at $x$} and we denote the linear map so defined by $f'(x)$. 

In this paragraph we consider the particular case when  $X=\R$ and we introduce further notions. In this setting, one can check that $f$ is \emph{semi-differentiable at $x\in D$} if and only if the left and right derivatives at $x\in D$ defined by
\Eq{*}{
  f'_-(x):=-f'(x;-1)=\lim_{t\to0-}\frac{f(x+t)-f(x)}{t}
  \quad\text{and} \quad
  f'_+(x):=f'(x;1)=\lim_{t\to0+}\frac{f(x+t)-f(x)}{t}  
}
exist. Then, we have
\Eq{f'}{
  f'(x;u)=f'_+(x)u \quad(u\geq0)
  \qquad\text{and} \qquad
  f'(x;u)=f'_-(x)u \quad(u\leq 0).
}
One can easily see that if $f$ is semi-differentiable at $x$, then $f$ is continuous at $x$. Clearly, $f$ is differentiable at $x$ if and only if $f'_-(x)=f'_+(x)$. 
We say that $f$ is \emph{nearly differentiable on $D$} if it is semi differentiable at each $x\in D$ and the set of points $x\in D$, where $f'_-(x)\neq f'_+(x)$ is at most countable. We call $f$ \emph{continuously semi-differentiable on $D$} if it is semi-differentiable on $D$ and, for all $x\in D$,
\Eq{*}{
  f'_-(x)
  =\lim_{u\to x-}f'_-(u)
  =\lim_{u\to x-}f'_+(u)
  \qquad\text{and} \qquad
  f'_+(x)
  =\lim_{u\to x+}f'_-(u)
  =\lim_{u\to x+}f'_+(u).
}
Obviously, if $f$ is differentiable on $D$, then it is continuously semi-differentiable on $D$ if and only if $f'$ is continuous on $D$, i.e., $f$ is continuously differentiable on $D$.

\Lem{Lip}{Let $(Y,\|\cdot\|)$ be a real normed space and let $D\subseteq\R$ be nonempty, open and let $f:D\to Y$ be continuously semi-differentiable on $D$. Then $f$ is nearly differentiable and locally Lipschitzian on $D$, i.e., for all $p\in D$, there exist $r>0$ and $L\geq0$ such that $[p-r,p+r]\subseteq D$ and 
\Eq{E:Lip}{
\|f(x)-f(y)\|\leq L|x-y| \qquad(x,y\in[p-r,p+r])
}
holds.}

\begin{proof} Observe that, for all $\varepsilon>0$, the set $\{x\in D \colon |f'_-(x)-f'_+(x)|\ge \varepsilon\}$ is discrete and therefore countable. Consequently, the set  $\{x\in D \colon f'_-(x)\ne f'_+(x)\}$ is also countable and therefore, $f$ is nearly differentiable.

Next, we prove that $f'_-$ and $f'_+$ are locally bounded over $D$. To show this, let $p\in D$ and choose $r>0$ so that the closed interval $I:=[p-r,p+r]$ be contained in $D$. We are going to verify that $f'_-$ is bounded over $I$. Assume, to the contrary, that $f'_-$ is not bounded over this interval. Then, for every $k\in\N$, there exists a point $x_k\in I$ such that $\|f'_-(x_k)\|\geq k$. We can extract a monotone subsequence $(y_\ell)=(x_{k_\ell})$ of $(x_k)$. Then, due to the compactness of $I$, this subsequence converges to a limit point $y_0\in I$. 
We have that $\|f'_-(y_\ell)\|\geq k_\ell$, therefore, the sequence $(\|f'_-(y_\ell)\|)$ tends to $+\infty$. This shows that $y_\ell\neq y_0$ for all $\ell\in\N$. 
If $(y_\ell)$ is increasing, then $y_\ell<y_0$ for all $\ell\in\N$, thus, by the continuous semi-differentiability of $f$, we get
\Eq{*}{
  \lim_{\ell\to\infty}\|f'_-(y_\ell)\|=\|f'_-(y_0)\|.
}
Similarly, if $(y_\ell)$ is decreasing, then we can see that
\Eq{*}{
  \lim_{\ell\to\infty}\|f'_-(y_\ell)\|=\|f'_+(y_0)\|.
}
Thus, we have obtained a contradiction in both cases, which then proves that $f'_-$ is bounded on $I$ and, since $p\in D$ was arbitrary, $f'_-$ is locally bounded on $D$. 

The proof of the local boundedness of $f'_+$ is completely analogous.

To show that $f$ is locally Lipschitzian, let $p\in D$ and choose $r>0$ so that $I:=[p-r,p+r]\subseteq D$. Let $L$ denote a common upper bound for both one-sided derivatives of $f$ over $I$. We are going to prove that the inequality \eq{E:Lip} holds. 

To the contrary, assume that \eq{E:Lip} is not valid.
Then, there exist $x_1,y_1\in I$ with $x_1<y_1$ and $\varepsilon>0$ such that
\Eq{*}{
  \|f(y_1)-f(x_1)\|>(L+\varepsilon)(y_1-x_1).
}

Using induction, we construct sequences $(x_n)$ and $(y_n)$ that satisfy $[x_{n+1},y_{n+1}] \subseteq [x_n,y_n]$,
\Eq{220}{
  \|f(y_n)-f(x_n)\|>(L+\varepsilon)(y_n-x_n)\\
}
and
\Eq{224}{
  y_n-x_n=2^{-n+1} (y_1-x_1)
}
for all $n \in \N$.

Clearly \eq{220} and \eq{224} hold for $n=1$. Assume that \eq{220} and \eq{224} are valid for some $n \in \N$ and $x_n,y_n \in I$ . Then, using the triangle inequality, we can easily see that 
\Eq{*}{
  \|f(y_n)-f(\tfrac{x_n+y_n}2)\|>(L+\varepsilon)\tfrac{y_n-x_n}2 \qquad\text{ or }\qquad
    \|f(\tfrac{x_n+y_n}2)-f(x_n)\|
    >(L+\varepsilon)\tfrac{y_n-x_n}2.
}
(Indeed, if the left hand side in both inequalities were not bigger than $(L+\varepsilon)\tfrac{y_n-x_n}2$,
then, by the triangle inequality, we would get
\Eq{*}{
  \|f(y_n)-f(x_n)\|
  \leq \|f(y_n)-f(\tfrac{x_n+y_n}2)\|
  +\|f(\tfrac{x_n+y_n}2)-f(x_n)\|
  \leq 2(L+\varepsilon)\tfrac{y_n-x_n}2,
}
which would contradict \eq{220}.)

Now we set $(x_{n+1},y_{n+1}):=(\tfrac{x_n+y_n}2,y_n)$ if the first inequality is valid and $(x_{n+1},y_{n+1}):=(x_n,\tfrac{x_n+y_n}2)$ otherwise. Then we see that \eq{220} holds with $n$ replaced by $n+1$. Furthermore, we have $y_{n+1}-x_{n+1}=\frac12 (y_n-x_n)=2^{-n}(y_1-x_1)$, which shows that \eq{224} is also valid with $n$ replaced by $n+1$. The inclusion $[x_{n+1},y_{n+1}] \subseteq [x_n,y_n]$ is also obvious by the construction.

In view of the inclusion and \eq{224}, we can see that $(x_n)$ is increasing and $(y_n)$ is decreasing, both sequences $(x_n)$ and $(y_n)$ converge to a common limit~$z \in I$. 

Let us now consider three cases. If $x_{n_0} = z$ for some $n_0 \in \N$, then we get $x_{n}=z<y_n$ for all $n \ge n_0$ and from the inequality \eq{220}, for all $n\geq n_0$, we obtain
\Eq{*}{
  \left\|\frac{f(y_n)-f(z)}{y_n-z}\right\|>L+\varepsilon.
}
Upon taking the limit $n\to\infty$,
we conclude that $\|f'_+(z)\|\ge L+\varepsilon$, a contradiction. Similarly, if $y_{n_0} = z$ for some $n_0 \in \N$, then we get $\|f'_-(z)\|\ge L+\varepsilon$, which is again a contradiction. 

In the final case, that is when $x_n< z<y_n$ for all $n\in \N$, by \eq{220} we obtain
\Eq{*}{
\|f(y_n)-f(z)\|+\|f(z)-f(x_n)\|&\ge \|f(y_n)-f(x_n)\|>(L+\varepsilon)(y_n-x_n)\\
&=(L+\varepsilon)((y_n-z)+(z-x_n)) \qquad (n \in \N)
}
Thus
$\|f(y_n)-f(z)\|\ge (L+\varepsilon)(y_n-z)$ for infinitely many $n \in\N$ or $\|f(z)-f(x_n)\|\ge (L+\varepsilon)(z-x_n)$ for infinitely many $n \in\N$. After taking a subsequence we reduce this case to one of the previous ones.
\end{proof}

The following two lemmas offer characterizations of monotonicity and convexity in terms of one-sided derivatives. 

\Lem{mon}{Let $I$ be a nonempty open real interval, let $f:I\to\R$ and assume that $f$ is continuous and differentiable from the left (resp.\ from the right) on $I$. Then $f$ is increasing if and only if $f'_-$ (resp. $f'_+$) is nonnegative on $I$.
}

For the proof we refer the reader to \cite[p.\ 13]{Bou04}.

\Lem{con}{Let $I$ be a nonempty open real interval, let $f:I\to\R$. Then $f$ is convex on $I$ if and only if it is continuously semi-differentiable on $I$ and the left (or, equivalently, the right) derivative of $f$ is an increasing function.}

For the proof we refer the reader to \cite[pp.\ 29--30]{Bou04} or to \cite[pp.\ 168--174]{Kuc85}.

In a subsequent result we are going to characterize the convexity of vector variable functions. 

\Lem{ChC}{
Let $(X,\|\cdot\|)$ be a real normed space, $D\subseteq X$ be a nonempty open convex set and $f:D\to\R$ be a continuously differentiable function and assume that $f':D\to X^*$ is semi-differentiable on $D$. Then $f$ is convex if and only if, for all $x\in D$, $u\in X$,
\Eq{E:fxuu}{
\big\langle(f')'(x;u),u\big\rangle\geq 0.
}
}

\begin{proof} The function $f$ is convex if and only if, for all $p,q\in D$, the real variable function 
\Eq{*}{
\varphi(s):=f((1-s)p+sq)=f(p+s(q-p))\qquad (s\in (0,1))
}
is convex. Due to the continuous differentiability of $f$, we have that $\varphi$ is also continuously differentiable on $(0,1)$ and we have that
\Eq{*}{
  \varphi'(s)=\big\langle f'(p+s(q-p)),q-p\big\rangle\qquad (s\in (0,1)).
}
Observe that
\Eq{*}{
  (\varphi')'_+(s)
  &=\lim_{t\to0+}
    \frac{\varphi'(s+t)-\varphi'(s)}{t}\\
  &=\lim_{t\to0+}
    \bigg\langle \frac{f'((p+s(q-p))+t(q-p))-f'(p+s(q-p))}{t},q-p\bigg\rangle\\
  &=\bigg\langle \lim_{t\to0+}\frac{f'((p+s(q-p))+t(q-p))-f'(p+s(q-p))}{t},q-p\bigg\rangle\\
  &=\big\langle (f')'(p+s(q-p);q-p),q-p\big\rangle.
}
The function $\varphi$ is convex if and only if $\varphi'$ is increasing, which, according to \lem{mon}, holds if and only if $(\varphi')'_+$ is nonnegative on $(0,1)$, i.e., for all $p,q\in D$, $s\in(0,1)$,
\Eq{*}{
  \big\langle (f')'(p+s(q-p);q-p),q-p\big\rangle\geq0.
}
The map $u\to \langle (f')'(p+s(q-p);u),u\rangle $
is quadratically positive homogeneous, therefore, the above inequality is equivalent to the validity of \eq{E:fxuu} for all $x\in D$ and $u\in X$.
\end{proof}

Finally, we deduce a chain rule in terms of semi-differentiability introduced above.

\Lem{CC}{
Let $(X,\|\cdot\|)$, $(Y,\|\cdot\|)$ and $(Z,\|\cdot\|)$ be real normed spaces, $D\subseteq X$ and $E\subseteq Y$ be nonempty open sets. Let $f:D\to E$ be semi-differentiable at $x\in D$ and let $g:E\to Z$ be locally Lipschitzian and semi-differentiable at $f(x)$. Then $g\circ f$ is semi-differentiable at $x\in D$ and
\Eq{E:CR}{
  (g\circ f)'(x;u)=g'(f(x);f'(x;u)) \qquad(u\in X).
}
}

\begin{proof} Using the local Lipschitz property of $g$ at $f(x)$, we can find $r>0$ and $L\geq0$ so that
\Eq{*}{
  \|g(y)-g(z)\|\leq L\|y-z\| \qquad(y,z\in B(f(x),r):=\{w\in Y\colon \|w-f(x)\|\leq r\}).
}
Let $u\in X$ be fixed. Then, we can choose $\tau>0$ so that
\Eq{*}{
  f(x+tu),f(x)+tf'(x;u)\in B(f(x),r) \qquad(t\in(0,\tau)).
}
Therefore,
\Eq{*}{
  \limsup_{t\to0+}\bigg\|\frac{g(f(x+tu))-g(f(x)+tf'(x;u))}{t}\bigg\|
  &= \limsup_{\tau>t\to0+}\bigg\|\frac{g(f(x+tu))-g(f(x)+tf'(x;u))}{t}\bigg\|\\
  &\le \limsup_{\tau>t\to0+} L\bigg\|\frac{f(x+tu)-f(x)}{t}-f'(x;u)\bigg\|=0,
}
which yields that
\Eq{*}{
  \lim_{t\to0+}\frac{g(f(x+tu))-g(f(x)+tf'(x;u))}{t}=0.
}
Consequently,
\Eq{*}{
  (g\circ f)'(x;u)
  &=\lim_{t\to0+}\frac{g(f(x+tu))-g(f(x))}{t}\\
  &=\lim_{t\to0+}\bigg(\frac{g(f(x+tu))-g(f(x)+tf'(x;u))}{t}+\frac{g(f(x)+tf'(x;u))-g(f(x))}{t}\bigg)\\
  &=\lim_{t\to0+}\frac{g(f(x)+tf'(x;u))-g(f(x))}{t}
  =g'(f(x);f'(x;u)).
}
This completes the proof of the chain rule stated in \eq{E:CR}.
\end{proof}



\section{\label{sec:aa} A Necessary Condition for the Subadditivity of General Means}

If $\M^{[n]}\colon I^n\to I$ is an $n$-variable mean and $\sigma$ is a permutation of the set $\{1,\dots,n\}$, then the mean $\M^{[n]}_\sigma\colon I^n\to I$ is defined by
\Eq{*}{
\M^{[n]}_{\sigma}(x_1,\dots,x_n):= \M^{[n]}(x_{\sigma(1)},\dots,x_{\sigma(n)})
\quad (x_1,\dots,x_n \in I).
} 
We say that $\M^{[n]}\colon I^n\to I$ is \emph{cyclically symmetric} if $\M^{[n]}=\M^{[n]}_\sigma$ for all cyclic permutations $\sigma$ of the index set $\{1,\dots,n\}$. 

\Thm{M-A}{
Let $n\geq 2$ and let $\M^{[n]}:\R_+^n\to\R_+$ be a cyclically symmetric subadditive (superadditive) $n$ variable mean. Then $\A^{[n]} \le \M^{[n]}$ (resp. $\M^{[n]} \le \A^{[n]}$). 
}

\begin{proof} 
Assume that $\M^{[n]}$ is subadditive and let $(x_1,\dots,x_n)\in\R_+^n$.
Using that $\M^{[n]}$ is cyclically symmetric, and then the subadditivity of $\M$, we obtain
\Eq{*}{
n\M^{[n]}(x_1,\dots,x_n)
&=\sum_{\sigma \in \Cyc_n} \M^{[n]}_\sigma(x_1,\dots,x_n)
=\sum_{\sigma \in \Cyc_n} \M^{[n]}(x_{\sigma(1)},\dots,x_{\sigma(n)})\\
&\ge \M^{[n]}\bigg(\sum_{\sigma \in \Cyc_n} x_{\sigma(1)},\dots,\sum_{\sigma \in \Cyc_n} x_{\sigma(n)}\bigg)\\
&= \M^{[n]}(x_1+\dots+x_n,\dots,x_1+\dots+x_n)=x_1+\dots+x_n.
}
This implies $\A^{[n]}(x)\le\M^{[n]}(x)$. As $(x_1,\dots,x_n)\in\R_+^n$ was taken arbitrarily, we can conclude that $\A^{[n]}\le\M^{[n]}$ holds.
\end{proof}

\section{Subadditivity of Quasi-arithmetic Means\label{sec:qa}}

In the next result we establish a complete characterization of the subadditivity of quasi-arithmetic means. To simplify the notation, given a differentiable function $f:D\to\R$ whose derivative is semi-differentiable, we shall write $f''_-$ and $f''_+$ for $(f')'_-$ and $(f')'_+$, respectively.
 
 \Thm{QAconvex}{
 Let $f \colon \R_+ \to \R$ be a continuous, strictly increasing function. Then the following assertions are equivalent to each other.
 \begin{enumerate}[(i)]
 \item For all $n\in\N$, the $n$-variable mean $\QA{f}^{[n]}$ is subaddditive;
 \item The $2$-variable mean $\QA{f}^{[2]}$ is subaddditive;
 \item The map $\Phi_f\colon f(\R_+)^2\to\R$ defined by
 \Eq{*}{
   \Phi_f(u,v):=f\big(f^{-1}(u)+f^{-1}(v)\big)
 }
 is concave.
 \item The function $f$ is differentiable with a positive derivative and the map $\Psi_f\colon \R_+^2\to\R$ defined by
 \Eq{*}{
   \Psi_f(x,y):=\frac{f(x)-f(y)}{f'(y)}
 }
 is subadditive;
 \item The function $f$ is differentiable with a positive derivative, its derivative is continuously semi-differentiable, and there exists an $\alpha\in[0,\infty]$ such that $f''_+$ is positive on $(0,\alpha)$ and $f''_+=0$ on $[\alpha,\infty)$, furthermore, $\frac{f'}{f''_+}$ is increasing and superadditive on $(0,\alpha)$, i.e.,
 \Eq{SA}{
  \frac{f'}{f''_+}(x+y)\geq \frac{f'}{f''_+}(x)+\frac{f'}{f''_+}(y) \qquad(x,y>0,\,x+y<\alpha).
 }
 \end{enumerate}
}
 
\begin{proof} The equivalence of the first three assertion is well-known from the theory of quasi-arithmetic means. The implication $(i)\Rightarrow(ii)$ is trivial. Assuming $(ii)$ and doing some simple substitutions, one could derive the Jensen concavity of $\Phi_f$. Using continuity, its concavity now follows. Finally, having the concavity of $\Phi_f$, one can obtain the $n$ variable Jensen inequality for $\Phi_f$, from which, after some obvious substitutions, the first assertion can be deduced.

Now assume that any of the equivalent assertions $(i)$, $(ii)$, and $(iii)$ are valid. Applying \thm{M-A} to the two variable mean $\QA{f}^{[2]}$, it follows that $\A^{[2]}\leq \QA{f}^{[2]}$. By the celebrated comparison theorem of quasi-arithmetic means, this inequality implies that $f$ is convex. Consequently, $f$ is continuously semi-differentiable, and the functions $f'_-$ and $f'_+$ are positive and increasing. 

Applying some general results concerning quasideviation means (cf. \cite[Theorem~11]{Pal88a}), it can be shown that $\QA{f}$ is subadditive if and only if there exist two functions $a,\,b \colon \R_+^2 \to \R$ such that
\Eq{E:condC}{
f(x+y)-f(u+v)\le a(u,v) (f(x)-f(u)) + b(u,v) (f(y)-f(v))
}
is valid for all $x,y,u,v \in I$.

For $x>u$ with $y=v$ we obtain
\Eq{*}{
\frac{f(x+v)-f(u+v)}{x-u} \le 
a(u,v) \cdot \frac{f(x)-f(u)}{x-u}.
}
Upon taking the limit $x \searrow u$, it follows that  
\Eq{*}{
  f'_+(u+v) \le a(u,v) f'_+(u).
}
Therefore
\Eq{*}{
\frac{f'_+(u+v)}{f'_+(u)}\leq a(u,v).
}
Analogously, we obtain
\Eq{*}{
a(u,v) \le \frac{f'_-(u+v)}{f'_-(u)},
}
which implies the double inequality
\Eq{E:DI}{
\frac{f'_+(u+v)}{f'_+(u)}\le a(u,v) \le \frac{f'_-(u+v)}{f'_-(u)}.
}

Take $x\in \R_+$ arbitrarily. As $f$ is convex, we know that $f'_-(x) \le f'_+(x)$. Using that $f$ is nearly differentiable, we can take $u \in (0,x)$ such that $f$ is differentiable at $u$. Then \eq{E:DI} with $v:=x-u$ simplifies to 
\Eq{*}{
\frac{f'_+(x)}{f'(u)} \le \frac{f'_-(x)}{f'(u)}
}
which implies $f'_+(x) \le f'_-(x)$. Consequently, $f$ is differentiable at $x$. As $f$ is convex, we get $f \in \mathscr{C}^1(I)$ with $f'>0$. Then, in view of \eq{E:DI} and the similar inequality for the function $b$, one can conclude that
\Eq{*}{
a(u,v) =\frac{f'(u+v)}{f'(u)}, \qquad b(u,v)=\frac{f'(u+v)}{f'(v)} \qquad (u,\, v \in \R_+).
}
Now condition \eq{E:condC} can be equivalently rewritten as 
\Eq{E:condC1}{
\frac{f(x+y)-f(u+v)}{f'(u+v)}\le \frac{f(x)-f(u)}{f'(u)} + \frac{f(y)-f(v)}{f'(v)}\qquad (x,y,u,v \in \R_+).
}
This implies that the function $\Psi_f$ is subadditive on $\R_+^2$ and completes the proof of assertion $(iv)$.

Now assume that assertion $(iv)$ holds, that is, $f$ is differentiable with a positive derivative and satisfies the inequality \eq{E:condC1}.
Therefore, for all $x,y,u,v \in \R_+$,
\Eq{*}{
f(x+y)-f(u+v)\le \frac{f'(u+v)}{f'(u)}(f(x)-f(u)) + \frac{f'(u+v)}{f'(v)}(f(y)-f(v)).
}
Interchanging the role of $(x,y)$ with $(u,v)$, we get
\Eq{*}{
f(u+v)-f(x+y)\le \frac{f'(x+y)}{f'(x)}(f(u)-f(x)) + \frac{f'(x+y)}{f'(y)}(f(v)-f(y)).
}
Adding up the above inequalities side by side, we can conclude that
\Eq{E:MD}{
0\le \bigg(\frac{f'(u+v)}{f'(u)}-\frac{f'(x+y)}{f'(x)}\bigg)(f(x)-f(u))+\bigg(\frac{f'(u+v)}{f'(v)}-\frac{f'(x+y)}{f'(y)}\bigg)(f(y)-f(v)).
}

Now substitute $v:=y$ into \eq{E:MD}. Then we obtain
\Eq{*}{
0\le \bigg(\frac{f'(u+y)}{f'(u)}-\frac{f'(x+y)}{f'(x)}\bigg)(f(x)-f(u))\qquad(x,y,u\in\R_+).
}
Due to the strict increasingness of $f$, this inequality shows that, for all $y\in\R_+$,
the map
\Eq{*}{
  \R_+\ni x\mapsto \frac{f'(x+y)}{f'(x)}
}
is decreasing. Therefore, for all $x,y\in\R_+$, we get
\Eq{*}{
  \frac{f'(x+y)}{f'(x)}\geq \frac{f'((x+y)+y)}{f'(x+y)}
}
Let $x<u$ be arbitrary. Then, with $y:=\frac{u-x}{2}$, the above inequality yields that
\Eq{*}{
  \Big(f'\Big(\frac{x+u}{2}\Big)\Big)^2
  \geq f'(x)f'(u) \qquad(x,u\in\R_+).
}
This inequality shows that function $g:=\log\circ f'$ is Jensen concave. Since $f'$ is Borel measurable, hence $g$ is also Borel measurable. According to the extension of the Bernstein--Doetsch theorem (cf.\ \cite{BerDoe15}) by Sierpiński \cite{Sie20}, we can conclude that $g$ is concave. Therefore, $g$ is continuously semi-differentiable, and the left and right derivatives of $g$ are decreasing, furthermore, we have 
\Eq{E:fg}{
\frac{f_-''}{f'}=g'_-\geq g'_+=\frac{f_+''}{f'}.
}
Consequently, $f'=\exp\circ g$ is also continuously semi-differentiable.
Since $f$ is convex, we have that $f'$ is increasing, therefore, $f''_-\geq f''_+\geq0$, which imply that $g'_-\geq g'_+\geq0$. 

Let $\alpha\in[0,\infty]$ be chosen so that $g'_+$ is positive on $(0,\alpha)$ and $g'_+=0$ on $(\alpha,\infty)$.
Therefore, by \eq{E:fg}, $f''_+$ is positive on $(0,\alpha)$ and $f''_+=0$ on $(\alpha,\infty)$. Due to the continuous semi-differentiability of $f$, we get that $f''_+(\alpha)=0$ is also valid. Using that the left and right derivatives of $g$ are decreasing, it follows from \eq{E:fg} that the functions $\frac{f_-''}{f'}$ and $\frac{f_+''}{f'}$ are decreasing, consequently
the functions $\frac{f'}{f_-''}$ and $\frac{f'}{f_+''}$
are increasing on $(0,\alpha)$.

Now we return to inequality \eq{E:MD}. Let $(x,y),(w,z)\in\R_+^2$ be fixed and substitute $u:=x+tw$ and $v:=y+tz$ into \eq{E:MD}, where $t>0$. Then, we get
\Eq{E:MD-t}{
0\ge &\bigg(\frac{f'(x+y+t(u+v))}{f'(x+tw)}-\frac{f'(x+y)}{f'(x)}\bigg)(f(x+tw)-f(x))
\\&+\bigg(\frac{f'(x+y+t(u+v))}{f'(y+tz)}-\frac{f'(x+y)}{f'(y)}\bigg)(f(y+tz)-f(y)).
}
By the differentiability of $f$, we have that
\Eq{*}{
  \lim_{t\to0+}\frac{f(x+tw)-f(x)}{t}=f'(x)w,\qquad
  \lim_{t\to0+}\frac{f(y+tz)-f(y)}{t}=f'(y)z.
}
On the other hand, using that $f'$ is right-differentiable everywhere, we can see that the maps $t\mapsto f'(x+tw)$ and $t\mapsto f'(x+y+t(u+v))$ are right differentiable at $t=0$, whence it follows that
\Eq{*}{
  \lim_{t\to0+}\frac{1}{t} \bigg(\frac{f'(x+y+t(u+v))}{f'(x+tw)}-\frac{f'(x+y)}{f'(x)}\bigg)
  =\frac{f''_+(x+y)(u+v)\cdot f'(x)-f'(x+y)\cdot f''_+(x)w}{f'(x)^2}.
}
Similarly, we can obtain that
\Eq{*}{
  \lim_{t\to0+}\frac{1}{t} \bigg(\frac{f'(x+y+t(u+v))}{f'(y+tz)}-\frac{f'(x+y)}{f'(y)}\bigg)
  =\frac{f''_+(x+y)(u+v)\cdot f'(y)-f'(x+y)\cdot f''_+(y)z}{f'(y)^2}.
}
Now dividing the inequality \eq{E:MD-t} by $t^2$ and then taking the limit as $t\to0+$, we can conclude that
\Eq{*}{
0\ge &\frac{f''_+(x+y)(u+v)\cdot f'(x)-f'(x+y)\cdot f''_+(x)w}{f'(x)^2}f'(x)w
\\&+\frac{f''_+(x+y)(u+v)\cdot f'(y)-f'(x+y)\cdot f''_+(y)z}{f'(y)^2}f'(y)z.
}
Dividing this inequality by $f'(x+y)$ side by side, it simplifies to
\Eq{*}{
   0\geq \frac{f''_+(x+y)}{f'(x+y)}(u+v)^2
     -\frac{f''_+(x)}{f'(x)}w^2
     -\frac{f''_+(y)}{f'(y)}z^2.
}
Assume now that $x,y>0$ with $x+y<\alpha$. Then $f''_+/f'$ is positive at $x,y$ and $x+y$. Define $(w,z)=\big(\frac{f''_+(y)}{f'(y)},\frac{f''_+(x)}{f'(x)}\big)$. Then the above inequality yields that
\Eq{*}{
   \frac{f''_+(x+y)}{f'(x+y)}\bigg(\frac{f''_+(x)}{f'(x)}+\frac{f''_+(y)}{f'(y)}\bigg)^2
     \leq\frac{f''_+(x)}{f'(x)}\bigg(\frac{f''_+(y)}{f'(y)}\bigg)^2
     +\frac{f''_+(y)}{f'(y)}\bigg(\frac{f''_+(x)}{f'(x)}\bigg)^2
}
which simplifies to
\Eq{*}{
   \frac{f''_+(x+y)}{f'(x+y)}\bigg(\frac{f''_+(x)}{f'(x)}+\frac{f''_+(y)}{f'(y)}\bigg)
     \leq\frac{f''_+(x)}{f'(x)}\frac{f''_+(y)}{f'(y)}.
}
Hence,
\Eq{*}{
   \frac{f'(x)}{f''_+(x)}+\frac{f'(y)}{f''_+(y)}
     \leq\frac{f'(x+y)}{f''_+(x+y)},
}
which shows that $f'/f''_+$ is superadditive over $(0,\alpha)$. This completes the proof of the implication $(iv)\Rightarrow(v)$.

In the last step of the proof, we show that assertion $(v)$ implies $(iii)$. The function $\Phi_f$ is continuously differentiable, therefore,
we can apply \lem{ChC} to verify its concavity. For this, we have to evaluate the directional derivative of $\Phi'_f$ first.

An elementary computation yields, for all $(X,Y)\in f(\R_+)^2$, that
\Eq{*}{
  \partial_1\Phi_f(X,Y)=\frac{f'(f^{-1}(X)+f^{-1}(Y))}{f'(f^{-1}(X))},\qquad 
  \partial_2\Phi_f(X,Y)=\frac{f'(f^{-1}(X)+f^{-1}(Y))}{f'(f^{-1}(Y))},
}
Since $f'$ is locally Lipschitzian, we can apply
the chain rule established in \lem{CC}. Then, for all $(X,Y)\in f(\R_+)^2$ and $(U,V)\in\R^2$, we get
\Eq{*}{
  (\partial_1\Phi_f)'((X,Y),(U,V))
  &=\frac{\big[(f')'\big(f^{-1}(X)+f^{-1}(Y);\frac{U}{f'(f^{-1}(X))}+\frac{V}{f'(f^{-1}(Y))}\big)\big]\cdot f'(f^{-1}(X))}{f'(f^{-1}(X))^2}\\
  &\qquad-\frac{f'(f^{-1}(X)+f^{-1}(Y))\cdot (f')'\big(f^{-1}(X);\frac{U}{f'(f^{-1}(X))}\big)}{f'(f^{-1}(X))^2}.
}
Thus, with the substitutions $(X,Y):=(f(x),f(y))$ and $(U,V):=(uf'(x),vf'(y))$, where $x,y\in\R_+$ and $u,v\in\R$, we get
\Eq{*}{
  (\partial_1\Phi_f)'((f(x),f(y)),(uf'(x),vf'(y)))
  &=\frac{(f')'(x+y;u+v)\cdot f'(x)-f'(x+y)\cdot (f')'(x;u)}{f'(x)^2}.
}
A similar argument yields that
\Eq{*}{
  (\partial_2\Phi_f)'((f(x),f(y)),(uf'(x),vf'(y)))
  &=\frac{(f')'(x+y;u+v)\cdot f'(y)-f'(x+y)\cdot (f')'(y;v)}{f'(y)^2}.
}

In view of \lem{ChC}, to prove that $\Phi_f$ is concave on $f(\R_+)^2$, we have to show that, for all $x,y\in\R_+$ and $u,v\in\R$,
\Eq{*}{
 \left\langle (\Phi_f')'\big((f(x),f(y));(uf'(x),vf'(y))\big),(uf'(x),vf'(y))\right\rangle\leq0.
}
Equivalently, we need to show that, for all $x,y\in\R_+$ and $u,v\in\R$,
\Eq{*}{
  (\partial_1\Phi_f)'((f(x),f(y)),(uf'(x),vf'(y)))\cdot uf'(x)+(\partial_2\Phi_f)'((f(x),f(y)),(uf'(x),vf'(y)))\cdot vf'(y)\leq0.
}
Using what we have established, we need to check that
\Eq{*}{
  &\frac{(f')'(x+y;u+v)\cdot f'(x)-f'(x+y)\cdot (f')'(x;u)}{f'(x)^2}\cdot uf'(x)\\
  &\qquad+\frac{(f')'(x+y;u+v)\cdot f'(y)-f'(x+y)\cdot (f')'(y;v)}{f'(y)^2}\cdot vf'(y)\leq0.
}
equivalently,
\Eq{539}{
  \frac{(f')'(x+y;u+v)}{f'(x+y)}\cdot(u+v)
  -\frac{(f')'(x;u)}{f'(x)}\cdot u
  -\frac{(f')'(y;v)}{f'(y)}\cdot v\leq0.
}

First we are going to prove the inequality
\Eq{546}{
 \frac{f''_+}{f'}(x+y)\cdot (u+v)^2
 \leq \frac{f''_+}{f'}(x)\cdot u^2+\frac{f''_+}{f'}(y)\cdot v^2
 \qquad (u,v\in\R,\,x,y\in\R_+).
}
If $x+y\geq\alpha$, then the left hand side is equal to $0$, therefore the inequality is trivial. 

Assume now that $x,y\in\R_+$ with $x+y<\alpha$. Then, by the superadditivity of $f'/f''_+$ over $(0,\alpha)$, we get
\Eq{*}{
\frac{f''_+}{f'}(x+y)
=\frac{1}{\frac{f'}{f''_+}(x+y)}
\leq\frac{1}{\frac{f'}{f''_+}(x)+\frac{f'}{f''_+}(y)}
=\frac{\frac{f''_+}{f'}(x)\cdot\frac{f''_+}{f'}(y)}{\frac{f''_+}{f'}(x)+\frac{f''_+}{f'}(y)}.
}
The obvious inequality 
\Eq{*}{
  \left(\frac{f''_+}{f'}(x)\cdot u-\frac{f''_+}{f'}(y)\cdot v\right)^2\geq0 \qquad(u,v\in\R,\,x,y\in\R_+,\ x+y<\alpha)
}
immediately implies that
\Eq{*}{
\frac{\frac{f''_+}{f'}(x)\cdot\frac{f''_+}{f'}(y)}{\frac{f''_+}{f'}(x)+\frac{f''_+}{f'}(y)}\cdot(u+v)^2
\leq \frac{f''_+}{f'}(x)\cdot u^2+\frac{f''_+}{f'}(y)\cdot v^2
\qquad(u,v\in\R).
}
This completes the proof of \eq{546}. 

The equalities in \eq{f'} applied to $f'$ instead of $f$ yield
\Eq{*}{
(f')'(z;w)=
\begin{cases}
 f''_+(z)\cdot w & \mbox{if } w\geq0,\\[2mm]
 f''_-(z)\cdot w & \mbox{if } w<0.
\end{cases}\qquad(z\in\R_+,\ w\in\R).
}
On the other hand, since $f'$ is continuously semi-differentiable, we can also get that
\Eq{zw}{
  \lim_{\varepsilon \to 0^+} f''_+(z+\varepsilon w)=
\begin{cases}
 f''_+(z) & \mbox{if } w\geq0,\\[2mm]
 f''_-(z) & \mbox{if } w<0
\end{cases}  \qquad(z\in\R_+,\ w\in\R).
}
Therefore,
\Eq{*}{
  (f')'(z;w)=\lim_{\varepsilon \to 0^+} f''_+(z+\varepsilon w)\cdot w
  \qquad(z\in\R_+,\ w\in\R).
}

Now let $x,y\in\R_+$, $u,v\in\R$ and choose $\varepsilon_0>0$ so that, for all $\varepsilon\in(0,\varepsilon_0)$, we have $x+\varepsilon u,y+\varepsilon v>0$. Then, applying \eq{546} with $x+\varepsilon u,y+\varepsilon v$ instead of $x,y$, we get
\Eq{*}{
 \frac{f''_+}{f'}(x+y+\varepsilon (u+v))\cdot (u+v)^2
 \leq \frac{f''_+}{f'}(x+\varepsilon u)\cdot u^2+\frac{f''_+}{f'}(y+\varepsilon v)\cdot v^2
 \qquad (\varepsilon\in(0,\varepsilon_0)).
}
Upon taking the limit $\varepsilon \to 0^+$ and using the continuity of $f'$ and the equality \eq{zw} for the pairs $(z,w)\in\{(x,u),(y,v),(x+y,u+v)\}$, we can conclude that 
\eq{539} holds. Therefore, the proof of the concavity of $\Phi_f$ and, consequently, the proof of the final implication $(v)\Rightarrow(iii)$ is complete.
\end{proof}



\end{document}